\documentclass[twoside,12pt]{article}

\usepackage{amssymb,amsmath}

\newenvironment{proof}{\begin{trivlist}\item[]{\it
Proof.}}{\hfill$\square$\end{trivlist}}

\newtheorem{theorem}{Theorem}

\newtheorem{definition}[theorem]{Definition}
\newtheorem{lemma}[theorem]{Lemma}

\def\mc{{\mathbb C}}

\def\omq{{\cal O}(M_q)}

\def\tr{{\rm Tr}}
\def\qx{{}_qX}

\newcommand\zum[2]{\sum_{\substack{#1\\#2}}}
\date{}
\begin{document}
\title{The traces of quantum powers commute}
\author{M Domokos
\footnote{This research was supported through a European Community
Marie Curie Fellowship. Partially supported by OTKA No. F 32325 and
T 34530. 
}\; 
and T H Lenagan
\footnote{Corresponding author. Fax: +44 131 650 6553; \ Tel: +44 131 650 5078;} 
\\ 
\\ 
R\'enyi Institute of Mathematics, Hungarian Academy of\\ 
Sciences, P.O. Box 127, 1364 Budapest, Hungary\\
E-mail: domokos@renyi.hu  
\\
\\
School of Mathematics, University of Edinburgh,\\
James Clerk Maxwell Building, King's Buildings,\\ 
Mayfield Road, Edinburgh EH9 3JZ, Scotland\\
E-mail: tom@maths.ed.ac.uk
} 
\maketitle

\begin{abstract}
The traces of the quantum powers of a generic quantum matrix pairwise 
commute. This was conjectured by Kaoru Ikeda, in connection with certain 
Hamiltonian systems. The proof involves Newton's formulae for quantum 
matrices, relating traces of quantum powers with sums of principal minors. 
\end{abstract} 



\bigskip 

The quasi-classical limit of the quantum group $GL_q(n)$ is the 
$n^2$-variable commutative polynomial algebra generated by 
$y_{ij}$, $1\leq i,j\leq n$, together with the Poisson bracket given by 
\[\{y_{ij},y_{kl}\}=(\theta(i,k)+\theta(j,l))y_{il}y_{kj},\]
where 
$\theta(i,j)=1$, if $i<j$; $0$, if $i=j$; $-1$, if $i>j$; 
see \cite{ku1}. 
This Poisson bracket was used in \cite{i1} to construct a Hamiltonian system, 
which can be reduced to the classical Toda lattice. 
The approach of \cite{i1} is based on the fact that 
\begin{equation}\label{eq:1}
\{\tr(Y^k),\tr(Y^m)\}=0\mbox{ \ for \ }k,m=1,2,\dots,
\end{equation} 
where $Y$ is the $n\times n$ matrix with $y_{ij}$ as the $(i,j)$ entry; 
thus $\tr(Y^k)$, $k=1,2,\dots$ can be considered as an involutive 
set of integrals. 
It is therefore natural to look for a quantum analogue of \eqref{eq:1}. 
This problem is addressed in \cite{ku2}, \cite{i2}, \cite{i3}. 
Recall that the coordinate ring $\omq$ 
of $n\times n$ quantum matrices is the $\mc$-algebra generated by 
$x_{ij}$, $1\leq i,j\leq n$, subject to the relations 
$x_{ij}x_{kl}-x_{kl}x_{ij}
=(q^{\theta(j,l)}-q^{-\theta(i,k)})x_{il}x_{kj}$, 
$1\leq i,j,k,l\leq n$, where $q$ is a fixed non-zero element of 
$\mc$ 
(in some papers $q=e^{-h/2}$ is considered to be an element of 
the ring of formal power series 
$\mc[[h]]$, and $\omq$ is the $\mc[[h]]$-algebra generated by the $x_{ij}$; 
this convention has to be used in the definition of the 
`quasi-classical limit'). 
The so-called $q$-multiplication of matrices was introduced in \cite{z}. 
If $A=(a_{ij})$, $B=(b_{ij})$ are $n\times n$ matrices with entries from 
$\omq$, then let $A\star B$ be the matrix whose $(i,j)$-entry is 
\[(A\star B)_{ij}=\sum_{k=1}^nq^{\theta(j,k)}a_{ik}b_{kj}.\]  
Note that $\star$ is not an associative multiplication. 
Write $X$ for the $n\times n$ matrix with $x_{ij}$ as the $(i,j)$ entry.   
Define $\qx^k$ for $k=0,1,2,\dots$ as follows. 
Set $\qx^0=I$, the identity matrix, $\qx^1=X$, 
$\qx^2=X\star X$, $\qx^3=X\star(X\star X)$, and recursively, 
$\qx^k=X\star(\qx^{k-1})$. 
Denote by $\tr(\qx^k)$ the usual trace (the sum of the diagonal entries) 
of $\qx^k$. When $q$ is specialized to $1$, the element $\tr(\qx^k)$ goes to 
$\tr(Y^k)$ from \eqref{eq:1}. 
Ikeda \cite{i2} has conjectured that the 
following  quantum analogue of \eqref{eq:1} holds: 
\begin{equation}\label{eq:2} 
[\tr(\qx^k),\tr(\qx^m)]=0 \mbox{ \ for \ }k,m=1,2,\dots, 
\end{equation} 
where $[a,b]=ab-ba$. The special cases of \eqref{eq:2} when 
$k=1$, $m=2$, and later, when $k=1$, $m$ arbitrary were verified in 
\cite{i2}, \cite{i3}. In the present paper we prove \eqref{eq:2} in general. 
Note that this explains the reason behind \eqref{eq:1}.   

We need to recall a related set of elements in $\omq$. 
The {\em quantum determinant} is 
\[\sigma_n:=\sum_{\pi \in S_n}(-q)^{l(\pi)}
x_{1,\pi(1)} \dots x_{n,\pi(n)},\]  
where $l(\pi)$ denotes the number of pairs $i<j$ with 
$\pi(i)>\pi(j)$ for a permutation $\pi$. 
Take $k$-element subsets $K,L$ of $\{1,\dots,n\}$. Then the $x_{ij}$ with 
$i\in K$, $j\in L$ generate a subalgebra of $\omq$,   
isomorphic to the coordinate ring of $k\times k$ quantum matrices. 
So we can form the quantum determinant of this subalgebra. The 
resulting element is denoted by $[K|L]$, and is called a 
$k\times k$ {\em quantum minor} in $\omq$. 
The sum $\sum_K [K|K]$ of principal $k\times k$ quantum minors is denoted by 
$\sigma_k$ for $k=1,\dots,n$. 
In particular, $\sigma_1=\tr(X)$. 
For the sake of notational convenience, 
we define $\sigma_{n+1},\sigma_{n+2},\dots$ to be zero. 

The elements $\sigma_1,\dots,\sigma_n$ pairwise commute by \cite{dl}. 
Denote by $R$ the subalgebra of $\omq$ generated by $\sigma_1,\dots,\sigma_n$. 
So $R$ is an $n$-variable commutative polynomial algebra. 
The commutativity of $R$ was obtained in \cite{dl} as a by-product 
of its Hopf algebraic interpretation: the $\sigma_k$ are the basic 
coinvariants for a coaction (a $q$-analogue of the adjoint action) 
of $GL_q(n)$ on $\omq$; in an alternative formulation, $R$ is the subset of 
cocommutative elements in the coalgebra $\omq$. 

To simplify notation, write $t_k=\tr(\qx^k)$, $k=1,2,\dots$. 
We shall show that the $t_k$ pairwise commute by proving that 
they are all contained in $R$. (Actually, $t_1,\dots,t_n$ turns out to be 
another generating set of the algebra $R$.) 
More precisely, we shall show that the sequences  
$\sigma_1,\sigma_2,\dots$ and $t_1,t_2,\dots$ are related by the same 
Newton's formulae as in the classical case $q=1$, 
see Lemma~\ref{lemma:newton}. 
This will be derived from the proof of a Cayley-Hamilton theorem due to 
Zhang \cite{z}. We present this proof in detail, because 
(P 1.2) on page 103 of \cite{z} is wrong, consequently, 
the formulae we need are not contained in \cite{z}. 

\begin{definition} We define $Z_k$, for $k = 0,1,2, \dots$, recursively, via
\begin{align*}
Z_0 &:= I\\
Z_k &:= X\star Z_{k-1} + (-1)^k I\sigma_k, 
\mbox{ \quad for }k\geq 1. 
\end{align*}
\end{definition} 

In order to prove the next lemma, we need to use a quantum Laplace
expansion in the following form for a quantum  minor $[K|L]$ with $i,r\in
K$. This relation can be easily obtained from 
Proposition 8 in 9.2.2 of \cite{ks}.
Set $\delta_{i,r}=1$ if $i=r$, and $0$, otherwise.  

\begin{equation}\label{eq:laplace} 
\delta_{i,r} [K|L] = \sum_{s\in L}\, (-q)^{l(s,L) -l(r,K)}
x_{is}[K\backslash r|L\backslash s], 
\end{equation} 
where $l(u,J)$ stands for the number of elements $j\in J$ with $u>j$, 
for an element $u$ and a subset $J$ of $\{1,\ldots,n\}$, 
and $K\backslash r$ is the difference of the sets $K$ and $\{r\}$.  
Later we shall abbreviate $l(u,\{j\})$ as $l(u,j)$, 
so $l(u,j)$ is $1$, if $u>j$, and $0$, otherwise.  

\begin{lemma}\label{lemma}
The $(i,j)$ entry of $Z_k$ is 
\[
(-1)^k \zum{|J| = k+1}{i,j\in J}
q^{\theta(i,j)}(-q)^{l(i,J) - l(j,J)}
[J\backslash j |J\backslash i]
\]
for $k = 0,1,\dots, n-1$, and $Z_n=Z_{n+1}=Z_{n+2}=\dots$ is the zero matrix. 
\end{lemma}

\begin{proof}
The proof is by induction on $k$, the starting case where $k=0$ is
easily checked. Take $ 0\leq k\leq n-2$ and assume that the formula holds
for $k$. First, we compute the $(i,j)$ entry of $X\star
Z_k$:  

\begin{eqnarray*}
(X\star Z_k)_{i,j} 
	&=& 
\sum_{s=1}^n\, x_{is}q^{\theta(j,s)}(Z_k)_{s,j}\\
	&=& 
\sum_{s=1}^n\, x_{is}q^{\theta(j,s)}(-1)^k\left\{ 
\zum{|J| = k+1}{s,j\in J}q^{\theta(s,j)}(-q)^{l(s,J) - l(j,J)}[J\backslash j
|J\backslash s]\right\}\\
	&=&
\zum{|J| = k+1}{j\in J}(-1)^k \left\{\sum_{s\in J}(-q)^{l(s,J) -
l(j,J)}x_{is}[J\backslash j|J\backslash s]\right\}\\
	&=&
\zum{|J| = k+1}{j\in J}(-1)^k  f(i,j,J)
\end{eqnarray*}
where we have set 
\[
f(i,j,J) = \sum_{s\in J}(-q)^{l(s,J) -l(j,J)}x_{is}
[J\backslash j|J\backslash s].
\]

If $i\in J$, then $f(i,j,J) = \delta_{i,j}[J|J]$ by \eqref{eq:laplace}; 
and so, 
\begin{eqnarray*}
(X\star Z_k)_{i,j} 
	&=& 
\zum{|J| = k+1}{i,j\in J}(-1)^k\delta_{i,j}[J|J] + 
\zum{|J| =k+1}{j\in J, i\not\in J}(-1)^kf(i,j,J).
\end{eqnarray*}

Hence, 
\begin{eqnarray*}
(Z_{k+1})_{j,j} 
	&=& 
\zum{|J| = k+1}{j\in J}(-1)^k[J|J] + (-1)^{k+1}\sum_{|J| = k+1}[J|J]\\ 
	&=& 
\zum{|J| = k+1}{j\not\in J}(-1)^{k+1}[J|J]\\
	&=&
\zum{|J| = k+2}{j\in J}(-1)^{k+1}[J\backslash j|J\backslash j]
\end{eqnarray*}
which agrees with the formula in the statement of the lemma.

If $i\neq j$, then
\begin{eqnarray*}
(Z_{k+1})_{i,j}
        &=& 
(X\star Z_k)_{i,j}\\
	&=&
\zum{|J| =k+1}{j\in J, i\not\in J}(-1)^kf(i,j,J)
\end{eqnarray*}

For $|J| = k+1$ with $j\in J$ and $i\not\in J$, set $K:= J\sqcup i$, 
the disjoint union of $J$ and $\{i\}$. 
By \eqref{eq:laplace}, 
\begin{eqnarray*} 
[K\backslash j|K\backslash i]
	&=&
\sum_{s\in K\backslash i}(-q)^{l(s,K\backslash i)-l(i, K\backslash j)}
x_{is} [K\backslash \{i,j\}|K\backslash \{i,s\}]\\
	&=&
\sum_{s\in J}(-q)^{l(s,J)-l(i, J\sqcup i\backslash j)}
x_{is} [J\backslash j|J\backslash s]\\
        &=&
(-q)^{l(j,J)-l(i, J\sqcup i\backslash j)}f(i,j,J).
\end{eqnarray*}
Hence, 
\[
f(i,j,J)= (-q)^{l(i, J\sqcup i\backslash j) -l(j, J)} [K\backslash j|
[K\backslash i],
\]
where $K= J\sqcup i$ when $j\in J$ and $i\not\in J$. 
Thus, for $i\neq j$, we see that 
\begin{eqnarray*}
(Z_{k+1})_{i,j} 
	&=&
\zum{|K| = k+2}{i,j\in K}(-1)^k(-q)^{l(i,K\backslash j)-l(j, K\backslash i)}
[K\backslash j|K\backslash i]\\
	&=&
(-q)^{l(j,i)-l(i,j)}\left\{
\zum{|J| = k+2}{i,j\in J}(-1)^k(-q)^{l(i,J) - l(j,J)}
[J\backslash j|J\backslash i]\right\},
\end{eqnarray*}
which is the formula in the lemma for $Z_{k+1}$, since 
$(-q)^{l(j,i)-l(i,j)} = -q^{\theta(i,j)}$, provided that $i\neq j$. 
So we have proved the lemma for $k=0,1,\dots,n-1$. 
Applying the case $k=n-1$ and \eqref{eq:laplace} once more, 
we obtain that $X\star Z_{n-1} = (-1)^{n-1}I\sigma_n$;
and so $Z_n = 0$. Consequently, we have $0=Z_n=Z_{n+1}=Z_{n+2}=\dots$. 
\end{proof}

From the definition of $Z_k$, it follows that 
\begin{equation} \label{eq:4} 
Z_k = \qx^k - \qx^{k-1}\sigma_1 + \qx^{k-2}\sigma_2 
+ \dots + (-1)^kI\sigma_k, 
\end{equation} 
where $\qx^{k-j}\sigma_j$ is the matrix obtained after multiplying each 
entry of $\qx^{k-j}$ by $\sigma_j$ from the right. 
Thus the assertion $Z_n=0$ in Lemma~\ref{lemma} 
is the (second) Cayley-Hamilton theorem in \cite{z}. 

By taking the trace of the equality \eqref{eq:4} 
we obtain
\[
\tr(Z_k) = t_k - t_{k-1}\sigma_1 + t_{k-2}\sigma_2 + \dots
+(-1)^{k-1}t_1\sigma_{k-1} + (-1)^kn\sigma_k.
\]
Now, for $k=0,1,\dots,n-1$ we have by Lemma~\ref{lemma} that 
\begin{eqnarray*}
(Z_k)_{i,i} 
	&=&
(-1)^k\zum{|K| = k}{i\not\in K} [K|K]
\end{eqnarray*}
and so $\tr(Z_k) = (-1)^k(n-k)\sigma_k$. For $k\geq n$ we have 
$\tr(Z_k)=\tr(0)=0$. Comparing the two expressions for
$\tr(Z_k)$ we see that Newton's formulae hold: 

\begin{lemma}\label{lemma:newton} For $k=1,2,\dots$ we have 
\[t_k-t_{k-1}\sigma_1 + t_{k-2}\sigma_2 + \dots +
(-1)^{k-1}t_1\sigma_{k-1} + (-1)^kk\sigma_k = 0.\] 
(Recall that $0=\sigma_{n+1}=\sigma_{n+2}=\dots$.) 
\end{lemma}

A result of \cite{dl} together with Lemma~\ref{lemma:newton} imply the 
following: 

\begin{theorem}\label{thm} 
We have $[t_k,t_m]=0$ for $k,m=1,2,\dots$. 
The elements $t_{n+1},t_{n+2},\dots$ can be expressed as 
polynomials of $t_1,\dots,t_n$. 
\end{theorem} 

\begin{proof} 
It follows from Lemma~\ref{lemma:newton} 
by induction on $k$ that $t_k$ is contained 
in $R=\mc[\sigma_1,\dots,\sigma_n]$ for all $k$. 
Since $R$ is commutative by \cite{dl}, the $t_k$ pairwise commute. 
Moreover, $R=\mc[t_1,\dots,t_n]$, since $\sigma_1,\ldots,\sigma_n$ 
can be expressed as a polynomial of $t_1,\ldots,t_n$. 
This implies the second statement. 
(The fact that $t_{n+1},t_{n+2},\dots$ are polynomials of 
$t_1,\dots,t_{n-1},\sigma_n$ was shown in \cite{i2} by an elaborate argument.) 
\end{proof} 


\end{document}